\documentclass[11pt]{article}
\usepackage{latexsym}

\usepackage{amsmath, amscd}  
\usepackage{amsfonts}
\usepackage{amssymb}

\newtheorem{thm}{Theorem}
\newtheorem{cor}[thm]{Corollary}
\newtheorem{lemma}[thm]{Lemma}

\newtheorem{prop}[thm]{Proposition}

\begin{document}
\small

\title{\bf 
       A convex body whose centroid and Santal\'o point are far apart
       \footnote{Keywords: Santal\'o point, . 2000 Mathematics Subject Classification: 52A20, 53A15 }}

\author{Mathieu Meyer, Carsten Sch\"utt and Elisabeth M. Werner 
\thanks{Partially supported by an NSF grant, a FRG-NSF grant and  a BSF grant}}

\date{}

\maketitle

\begin{abstract}
We give an example of a convex body whose centroid and 
Santal\'o point are ``far apart".

\end{abstract}

\newpage
\section{Introduction. The main theorem.}

In his survey paper \cite{Gruenbaum1963},  Gr\"unbaum introduced measures of asymmetry
for convex bodies in $\mathbb R^{n}$,  i.e. compact, convex sets in $\mathbb R^{n}$ with nonempty interior.
Measures of asymmetry  determine the degree of  non-symmetricity of a non symmetric  convex body.
Recall that a convex set $K$ in $\mathbb R^{n}$ is 
centrally symmetric with respect to the origin if $x\in K$ implies $-x\in K$.
More generally, $K$ is centrally symmetric with respect to the point
$z$ if $x\in K$ implies $2z-x\in K$.
\par
Gr\"unbaum argues that a measure of asymmetry should be a function
that is $0$ for centrally symmetric convex bodies and maximal for the
simplex and only for the simplex.
However, the simplex has enough symmetries to ensure that
all affine invariant points are identical. This leads to the question:
How far apart can specific affine invariant points be? 
This distance can be used as a measure for the asymmetry of the convex body.
\par
Of particular interest are the centroid $g(K)$ and the Santal\'o point  $s(K)$ of a convex body
$K$ in $\mathbb{R}^n$. These points are of major importance 
and have been widely studied. For the (less well known) Santal\'o point see e.g. 
\cite{Ball3,  Meyer-Pajor, MeyerWerner1998, Reisner, Schneider1993}.
Recent developments in this direction include e.g. the works 
\cite{ArtsteinKlartagMilman, Ball1, Ball2, BoeroeczkySchneider,FradeliziMeyer, Lehec, 
NazarovPetrovRyaboginZvavitch}.
We show here that the points can be very far apart, namely
\vskip 3mm

\begin{thm}\label{MainTheorem}
There is an absolute  constant $c>0$ and $n_{0}\in\mathbb N$ such that for all
$n\geq n_{0}$ there is a convex body $C=C_n$ in $\mathbb R^{n}$
with
\begin{equation}\label{MainTheorem1}
c
\leq\frac{\|\operatorname{g}(C)-\operatorname{s}(C)\|_{2}}
{\operatorname{vol}_{1}(\ell\cap C)}= \frac{\|\operatorname{g}(C)-\operatorname{s}(C)\|_{2}}
{w_C(u)} \ .
\end{equation}
$\|\cdot \|_{2}$ is the Euclidean norm, $l$ is the line through $g=g(C)$ and $s=s(C)$ and  
$w_C(u) = h_C(u) + h_C(-u)$ is the width of $C$ in directon of the unit vector $u=\frac{g-s}{\|g-s\|_{2}}$.
\end{thm}
\vskip 3mm      
\noindent
The proof actually shows that we can asymptotically determine the constant $c$ of the theorem.
$$
\frac{\|\operatorname{g}(C)-\operatorname{s}(C)\|_{2}}
{\operatorname{vol}_{1}(\ell\cap C)}
$$
is asymptotically with respect to the dimension,  greater than or equal to
$$
\left(1-\frac{1}{e}\right)\frac{\sqrt{e\pi}-2}{\sqrt{e\pi}+\frac{2}{e-1}}
= 0.142673...
$$
and is asymptotically with respect to the dimension,  less than or equal to
$$
\left(1-\frac{1}{e}\right)\frac{\sqrt{e}-1}{\sqrt{e}+\frac{1}{e-1}}
= 0.18383...
$$

\vskip 3mm

\vskip 3mm
\noindent
Throughout the paper we use the following notations.
$B_2^n(a,r)$ is the Euclidean ball in $\mathbb{R}^n$ centered at $a$ with
radius $r$. We denote $B_2^n=B_2^n(0,1)$. $S^{n-1}$ is the boundary
$\partial B_{2}^{n}$ of the Euclidean unit ball.
For $x \in \mathbb{R}^n$, let 
$\|x\|_p=\left(\sum _{i=1}^n |x_i|^p\right)^\frac{1}{p}$, 
if  $1 \leq p < \infty$ and $\|x\|_\infty= \max_{1 \leq i \leq n} |x_i|$, 
if $p=\infty$ and let $B_p^n=\{x \in \mathbb{R}^n: \|x\|_p \leq 1\}$  
be the unit ball of the space $l_p^n$.
\par
\noindent
For  a convex body $K$ in $\mathbb R ^n$ we denote the volume by 
$\mbox{vol}_n(K)$ (if we want to emphasize the dimension) or by $|K|$.  
$\mbox{co}[A,B]= \{\lambda a+(1-\lambda)b: a \in A, b\in B, 0 \leq \lambda \leq 1\}$ is the convex 
hull of $A$ and $B$.
\par 
$h_K(u)= \max_{x\in K} \langle u,x \rangle $ is the supportfunction of $K$ in direction $u$.
\par
For  $x\in K$, $K^x = (K-x)^\circ=\{y\in
{\mathbb R}^n: \langle y,z-x\rangle
\leq 1
\hspace{.1in}\mbox{for all z} \in K \}$ is the polar body of $K$ with
respect to $x$. 
\par
Let $\xi\in\mathbb R^{n}$ with $\|\xi\|_2=1$ and $s\in\mathbb R$.
Then the $(n-1)$-dimensional section of $K$ orthogonal to $\xi$ through $s\xi$ is
$$
K(s,\xi)
=\{x\in K| \langle \xi,x \rangle=s\}.
$$
We just write  $K(s)$ if it is clear  which direction $\xi$ is meant.
\vskip 3mm
\noindent
The
centroid $g(K)$  of a convex body
$K$ in
$\mathbb R^{n}$ is the point
$$
g(K)
=\frac{1}{\operatorname{vol}_{n}(K)}\int_{K}xdx.
$$
\par
\noindent
The Santal\'o point  $s(K)$ of a convex body $K$ is the unique point
$s(K)$
for which 
$$
\operatorname{vol}_{n}(K)\operatorname{vol}_{n}(K^{x})
$$
attains its minimum. It is more difficult to compute
the Santal\'o point as it is defined implicitly.
\par
Before we go into the construction of the convex bodies
that fulfill the properties of Theorem \ref{MainTheorem},
we discuss an example that shows that a more elaborate construction is necessary
to get the statement of Theorem \ref{MainTheorem}.
\par For instance, it is not enough to take one half of a Euclidean ball
$$
B=\{x\in\mathbb R^{n}|\|x\|_{2}\leq 1, x_{1}\geq 0\}\ .
$$
We compute the centroid of $B$. The coordinates
of the centroid are
$g(B)(2)=\cdots=g(B)(n)=0$ and 
\begin{eqnarray*}
g(B)(1)
&=&\frac{2}{|B_{2}^{n}|}
\int_{0}^{1}t\operatorname{vol}_{n-1}
\left(\sqrt{1-t^{2}}B_{2}^{n-1}\right)dt  
=\frac{2\operatorname{vol}_{n-1}  
\left(B_{2}^{n-1}\right)}{(n+1)|B_{2}^{n}|}\ .
\end{eqnarray*}
Now we compute the Santal\'o point of $B$.
The polar body of the Euclidean ball $B_{2}^{n}(\lambda e_{1},1)$
with center $\lambda e_{1}$, $0\leq\lambda<1$, and radius $1$ is
\begin{equation}\label{HalfBall1}
\left\{x\in \mathbb R^{n}\left|(1-\lambda^{2})^{2}
\left(x_{1}+\frac{\lambda}{1-\lambda^{2}}\right)^{2}+(1-\lambda^{2})x_{2}^{2}+\cdots+(1-\lambda^{2})x_{n}^{2}\leq 1\right.\right\}\ .
\end{equation}
Therefore, the polar body of the half ball $B$ with respect to $\lambda e_{1}$ is
the convex hull of the vector $-\frac{e_{1}}{\lambda}$ and the ellipsoid
(\ref{HalfBall1}).
We estimate the volume of $(B-\lambda e_{1})^{\circ}$.
Since $(B-\lambda e_{1})^{\circ}$ contains the ellipsoid
(\ref{HalfBall1})
$$
(1-\lambda^{2})^{-\frac{n+1}{2}}
\operatorname{vol}_{n}(B_{2}^{n})
\leq
\operatorname{vol}_{n}((B-\lambda e_{1})^{\circ}).
$$
On the other hand, as $(B-\lambda e_{1})^{\circ}$ is contained in the union
of the ellipsoid (\ref{HalfBall1}) and a cone with height $\frac{1}{\lambda}$
and with a base that is a Euclidean ball with radius $(1+\lambda^{2})^{-1}$,
$$
\operatorname{vol}_{n}((B-\lambda e_{1})^{\circ})
\leq\frac{\operatorname{vol}_{n}(B_{2}^{n})}{(1-\lambda^{2})^{\frac{n+1}{2}}}
+\frac{\operatorname{vol}_{n-1}(B_{2}^{n-1})}{n\lambda(1+\lambda^{2})^{n-1}}
$$
For $\lambda=\frac{1}{\sqrt{n}}$ we get
$$
\operatorname{vol}_{n}((B-\lambda e_{1})^{\circ})
\leq\left(1-\frac{1}{n}\right)^{-\frac{n+1}{2}}\operatorname{vol}_{n}(B_{2}^{n})
+\frac{\operatorname{vol}_{n-1}(B_{2}^{n-1})}{\sqrt{n}(1+\frac{1}{n})^{n-1}}\ .
$$
Since
$$
\frac{\operatorname{vol}_{n-1}(B_{2}^{n-1})}
{\operatorname{vol}_{n}(B_{2}^{n})}
\sim\sqrt{\frac{\pi n}{2}}
$$
we get
$$
\operatorname{vol}_{n}((B-\lambda e_{1})^{\circ})
\lesssim \left(\sqrt{e}
+\frac{1}{e}\sqrt{\frac{\pi}{2}}\right)\operatorname{vol}_{n}(B_{2}^{n})\ .
$$
On the other hand, for $\lambda=\frac{\gamma}{\sqrt{n}}$
$$
\operatorname{vol}_{n}((B-\lambda e_{1})^{\circ})
\geq\left(1-\frac{\gamma^{2}}{n}\right)^{-\frac{n+1}{2}}
\operatorname{vol}_{n}(B_{2}^{n}).
$$
With $1-t\leq e^{-t}$
$$
\operatorname{vol}_{n}((B-\lambda e_{1})^{\circ})
\geq e^{\gamma^{2}\frac{n+1}{2n}}
\operatorname{vol}_{n}(B_{2}^{n})
\ .
$$
Therefore, there is a $\gamma$ such that for all $n\in\mathbb N$
we have $s(B)(1)\leq\frac{\gamma}{\sqrt{n}}$. Thus 
$
\frac{\|\operatorname{g}(B)-\operatorname{s}(B)\|_{2}}
{\operatorname{vol}_{1}(\ell\cap B)}
$
is of order $\frac{1}{\sqrt{n}}$.
\vskip 3mm
The next lemma is well known (see \cite{Schneider1993}).
\vskip 3mm
\begin{lemma} \label{SantaloPoint}
For any convex body $K$, an interior point $x$ of $K$ is the
Santal\'o point if and only if $0$ is the centroid
of $(K-x)^{\circ}$.
\end{lemma}
\vskip 3mm

\noindent
This lemma can be rephrased as follows:
\vskip 2mm
{\em Let $K$ be a convex body. 
Then $0$ is the Santal\'o point of $K^{g(K)}$.}
\vskip 2mm
\noindent
Indeed, $K^{g(K)}=\left(K-g(K)\right)^{\circ}$ and 
$(K-g(K))^{\circ\circ}=K-g(K)$. 
Since $0$ is the centroid of $(K-g(K))^{\circ\circ}$,
it follows by Lemma \ref{SantaloPoint} that $0$ is the Santal\'o point
of $(K-g(K))^{\circ}=K^{g(K)}$.

\vskip 4mm
\noindent
For convex bodies
$K$ and
$L$ in
$\mathbb R ^n$ and natural numbers
$k$,
$0 \leq k \leq n$, the coefficients 
$$
V_{n-k,k}(K,L)=V( \underbrace {K,\dots,K}_{n-k} ,
\underbrace {L,\dots,L}_{k} )
$$ 
in the expansion
\begin{equation}\label{mixed volumes}
\mbox{vol}_{n}(\lambda _1 K+ \lambda _2 L)=\sum _{k=0}^n {{n} \choose
{k}}
\lambda _1^{n-k} \lambda _2^k V_{n-k,k}(K,L)
\end{equation}
are the mixed
volumes of $K$ and $L$ (see \cite{Schneider1993}).
\vskip 3mm
Now we introduce the convex bodies which will serve as candidates for Theorem \ref{MainTheorem}.
For convex bodies $K$ and $L$ in $\mathbb R ^n$ and real numbers $a >0$
and $b > 0$, we construct a convex body 
$M_n$ in $\mathbb R ^{n+1}$
\begin{equation}\label{definition:M}
M_n=\operatorname{co}[(K,-a),(L,b)]
=\{t(x,-a)+(1-t)(y,b)|x\in K,y\in L, 0\leq t\leq 1\}.
\end{equation}
The bodies we are using in Theorem \ref{MainTheorem} will be the polar bodies
to $M_{n}$.
The polar of $M_n$ can be described as follows: 
For  $-\frac{1}{a} \leq s \leq \frac{1}{b}$, the sections of $M_n^\circ$ orthogonal to $e_{n+1}$ 
and containing $se_{n+1}$ are
\begin{equation}\label{Mo(s)}
M_n^\circ(s) = M_n^\circ(s, e_{n+1}) =(1+s \ a)K^\circ \cap (1-s \ b)L^\circ.
\end{equation}
We show this.
\begin{eqnarray*}
M_n^{\circ}
&=&\{(z,s)\in\mathbb R^{n}\times \mathbb R|
\forall x\in K:<z,x>-sa\leq 1
\hskip 1mm\mbox{and}\hskip 1mm
\forall y\in L:<z,y>+sb\leq	1\}   \\
&=&\{(z,s)\in\mathbb R^{n}\times \mathbb R|
z\in (1+sa)K^{\circ}\cap(1-bs)L^{\circ}\}
\end{eqnarray*}
\vskip 4mm
The body $M_n$  is the convex hull of two 
$n$-dimensional faces $K$ and $L$.
In the following proposition we choose specific bodies for those faces.
We choose them   in such a way that
their volume product differs greatly, which has as effect that the 
centroid  and the Santal\'o point  of
$M_n^{g(M_n)}$ are ``far apart".  One face is  the
Euclidean ball and the other the unit ball of $\ell_{\infty}^{n}$.
We normalize both balls so that their volume is $1$.

\begin{prop}\label{thm:distance}
Let  
$$
K=\frac{B_2^n}{|B_2^n|^{\frac{1}{n}}}
\hskip 20mm
L=\frac{1}{2} \ B_\infty^{n}.
$$
and 
$a=1$ and $b=\frac{1}{e-1}$.
Let $M_n$ be the convex body in $\mathbb R ^{n+1}$ defined in
(\ref{definition:M}). Then
\vskip 3mm
(i) $\lim_{n \to \infty} \ g(M_n)=0$.
\vskip 3mm
(ii)
Let 
$$
s_{0}=-\frac{\sqrt{e}-1}{\sqrt{e}+\frac{1}{e-1}}
=-0.290815... 
\hskip 20mm
s_{1}
=
\frac{2-\sqrt{e\pi}}{\sqrt{e\pi}+\frac{2}{e-1}}
=-0.225705...
$$ 
Then for every $\varepsilon>0$ there is $n_{0}$ such that for all $n\geq n_{0}$ 
the $n+1$-st coordinate of the centroid
$g(M_{n}^{\circ})$ of
$M_{n}^{\circ}$ satisfies
$$
s_{0}-\varepsilon
\leq g(M^{\circ}_{n})(n+1)
\leq s_{1}+\varepsilon.
$$
The Santal\'o point of
$M_n^{g(M_n)}$, $s_n^*=s(M_n^{g(M_n)})=0$ 
 and
therefore the centroid $g_n^*$ and the Santal\'o point $s_n^*$ of
$M_n^{g(M_n)}$ satisfy
$$
\frac{\sqrt{e}-1}{\sqrt{e}+\frac{1}{e-1}}
\geq
\limsup_{n \to \infty} | g_n^* -s_n^* | 
\geq\liminf_{n \to \infty} | g_n^* -s_n^* | 
=\liminf_{n \to \infty} | g_n^* | 
\geq\frac{2-\sqrt{e\pi}}{\sqrt{e\pi}+\frac{2}{e-1}}\ .
$$
\end{prop}
\vskip 3mm
As for  the proof of Proposition \ref{thm:distance},   it follows from Lemma  \ref{SantaloPoint} that  $s(M_n^{g(M_n)})=0$. For the centroid $g_n^*$ 
we actually estimate the coordinates of $g(M_n^{\circ})$ and not
$g(M_n^{g(M_n)})$. Since $\lim g(M_n)=0$ it suffices to use a
perturbation argument. The remaining part of the proof of the proposition is in the next section.
\vskip 3mm
\noindent
{\bf Proof of
Theorem \ref{MainTheorem}.}  The proof follows immediately from Proposition \ref{thm:distance}: For every $n$, let $C_n=M_n^\circ$.
\vskip 5mm

\section{Proof of Proposition \ref{thm:distance}}
This section is devoted to the proof of Proposition \ref{thm:distance}. 
We will need several lemmas and then give the proof of the proposition at the end of the section.
\vskip 3mm

\begin{lemma}\label{lemma:gravityM}
Let $K$ and $L$ be convex bodies in $\mathbb R^n$ such that their
centroid are at $0$. Let $c > 0$ and let $M_n$ be the convex body in $\mathbb{R}^{n+1}$
$$
M_n=\mbox{co}[ (K,0),(L,c)].
$$
Then the $(n+1)$-coordinate of the centroid $g(M_n)(n+1)$ of $M_n$ satisfies
$$
g(M_n)(n+1)=\langle g(M_n),e_{n+1}\rangle =\frac{c}{n+2} \frac{\sum_{k=0}^n
(k+1)V_{n-k,k}(K,L)}{\sum_{k=0}^n V_{n-k,k}(K,L)}.
$$
\end{lemma}
\vskip 3mm

\noindent
{\bf Proof.}
By definition
$$
\langle g(M_n),e_{n+1} \rangle =\frac{\int_0^c w \  \mbox{vol}_n(M_n(w))dw}{\int_0^c
\mbox{vol}_n(M_n(w))dw}.
$$
Note that
$$
\mbox{co}[ (K,0),(L,c)]=\left\{\left(\left.(1-\frac{w}{c}) K +\frac{w}{c} L,w\right)   \right|0 \leq w
\leq c \right\}.
$$ 
Therefore 
\begin{eqnarray*}
\langle g(M_n),e_{n+1}\rangle&=&\frac{\int_0^c w \ \mbox{vol}_n((1-\frac{w}{c}) K
+\frac{w}{c} L)dw}{\int_0^c
\mbox{vol}_n((1-\frac{w}{c}) K +\frac{w}{c} L)dw}\\
&=&
c \ \frac{\int_0^1 t \  \mbox{vol}_n((1-t) K +t L)dt}{\int_0^1
\mbox{vol}_n((1-t) K +t L)dt}.
\end{eqnarray*}  
Now we use  the mixed volume formula (\ref {mixed volumes}). 
Thus
\begin{eqnarray*}
\langle g(M_n),e_{n+1}\rangle
&=&c \ \frac{\int_0^1 \sum_{k=0}^n  \biggl({{n} \choose
{k}} \  t^{k+1} \ (1-t)^{n-k} \ V_{n-k,k}(K,L) \biggr)\ dt}
{\int_0^1 \biggl(\sum_{k=0}^n  {{n} \choose
{k}} \  t^k \ (1-t)^{n-k} \ V_{n-k,k}(K,L)\biggr) \ dt} \\
&=&\frac{c}{n+2} \frac{\sum_{k=0}^n
(k+1)V_{n-k,k}(K,L)}{\sum_{k=0}^n V_{n-k,k}(K,L)},
\end{eqnarray*}
where we have also used the Betafunction
$$
B(k,l)= \int_0^1 t^ {k-1}(1-t)^{l-1}
dt=\frac{\Gamma (k)\Gamma (l)}{\Gamma (k+l)}, \ \ \  \ \ \  k,l > 0.
$$ 
$\square$
\vskip 3mm

\noindent
The next lemma is well known (\cite{Hadwiger}, p. 216, formula 54).
\vskip 3mm

\begin{lemma}\label{lemma:mixed volumes}
For all $n\in\mathbb N$ and $t\geq 0$
$$
\operatorname{vol}_n \left(B_2^n+t B_{\infty}^n \right) 
= \sum_{k=0}^n {n\choose k} 2^k  \operatorname{vol}_{n-k} \left(B_2^{n-k} \right) \ t^k,
$$
with the convention that  $\operatorname{vol}_{0}(B_{2}^{0})=1$.
Therefore, for $0 \leq k \leq n$,  
$$ 
V_{n-k,k}(B_2^n,B_{\infty}^n)=2^k\  \operatorname{vol}_{n-k} \left(B_2^{n-k} \right) .
$$
\end{lemma}
\vskip 3mm

\begin{lemma}\label{centroid}
The following formula holds
$$
\lim_{n\to\infty}
\frac{1}{n+2}\  \frac{\sum_{k=0}^n \frac{k+1}{\Gamma(1+\frac{n}{2})^{\frac{k}{n}}
\Gamma(1+\frac{n-k}{2})} }{\sum_{k=0}^n \frac{1}{\Gamma(1+\frac{n}{2})^{\frac{k}{n}}
\Gamma(1+\frac{n-k}{2})} }
=1-\frac{1}{e}.
$$
\end{lemma}

\vskip 3mm\noindent
{\bf Proof.}
One has 
$$\frac{1}{n+2}\  \frac{\sum_{k=0}^n \frac{k+1}{\Gamma(1+\frac{n}{2})^{\frac{k}{n}}
\Gamma(1+\frac{n-k}{2})} }{\sum_{k=0}^n \frac{1}{\Gamma(1+\frac{n}{2})^{\frac{k}{n}}
\Gamma(1+\frac{n-k}{2})} }= 
\frac{1}{n+2} \ \frac {\sum_{k=0}^n 
 \frac{n-k+1}{\Gamma(1+\frac{n}{2})^{\frac{n-k}{n}} \Gamma(1+\frac{k}{2})} }
{\sum_{k=0}^n 
  \frac{1}{\Gamma(1+\frac{n}{2})^{\frac{n-k}{n}} \Gamma(1+\frac{k}{2}) } }$$
  $$
=
\frac{1}{n+2} \frac{\sum_{k=0}^n 
\frac{  (n-k+1)\Gamma(1+\frac{n}{2})^{\frac{k}{n}} } {\Gamma(1+\frac{k}{2}) } }
{\sum_{k=0}^n 
\frac{\Gamma(1+\frac{n}{2})^{\frac{k}{n}} } {\Gamma(1+\frac{k}{2})} }
= \frac{n+1}{n+2} - \frac{1}{n+2} \frac{\sum_{k=0}^n 
k\frac{  \Gamma(1+\frac{n}{2})^{\frac{k}{n}} } {\Gamma(1+\frac{k}{2}) } }
{\sum_{k=0}^n 
\frac{\Gamma(1+\frac{n}{2})^{\frac{k}{n}} } {\Gamma(1+\frac{k}{2})} }.$$
It is thus needed to prove that 
\begin{equation}\label{centroid1}
\lim_{n\to\infty}
\frac{1}{n} \frac{\sum_{k=0}^n 
k\frac{  \Gamma(1+\frac{n}{2})^{\frac{k}{n}} } {\Gamma(1+\frac{k}{2}) } }
{\sum_{k=0}^n 
\frac{\Gamma(1+\frac{n}{2})^{\frac{k}{n}} } {\Gamma(1+\frac{k}{2})} }
= \frac{1}{e}.
\end{equation}
For every $x\geq0$, 
\begin{eqnarray*}
\sum_{k=0}^n \frac{kx^k}{\Gamma(1+\frac{k}{2})}
&=&\sum_{k=1}^n \frac{kx^k}{\frac{k}{2}\Gamma(\frac{k}{2})}
= 2x\sum_{k=1}^n \frac{x^{k-1}}{\Gamma(\frac{k}{2})}   \\
&= &2x\Big(\frac{1}{\Gamma(\frac{1}{2})}+\sum_{k=2}^{n} \frac{x^{k-1}}{\Gamma(\frac{k}{2})}\Big)
=2x\Big(\frac{1}{\Gamma(\frac{1}{2})}+x\sum_{k=0}^{n-2} \frac{x^k}{\Gamma(1+\frac{k}{2})}\Big).
\end{eqnarray*}
Thus for $x=\Gamma(1+\frac{n}{2})^{\frac{1}{n}}$
\begin{equation}\label{centroid2}
\sum_{k=0}^n \frac{k\Gamma(1+\frac{n}{2})^{\frac{k}{n}}}{\Gamma(1+\frac{k}{2})}
=2\Gamma\left(1+\frac{n}{2}\right)^{\frac{1}{n}}\left(\frac{1}{\Gamma(\frac{1}{2})}
+\Gamma\left(1+\frac{n}{2}\right)^{\frac{1}{n}}
\sum_{k=0}^{n-2} \frac{\Gamma(1+\frac{n}{2})^{\frac{k}{n}}}{\Gamma(1+\frac{k}{2})}\right).
\end{equation}
Let 
$$
A_n=\sum_{k=0}^n 
k\frac{  \Gamma(1+\frac{n}{2})^{\frac{k}{n}} } {\Gamma(1+\frac{k}{2}) }  \hbox{ and } 
B_n=\sum_{k=0}^n 
\frac{  \Gamma(1+\frac{n}{2})^{\frac{k}{n}} } {\Gamma(1+\frac{k}{2}) }.
$$
(\ref{centroid1}) is equivalent to
$$
\lim_{n\to\infty}\frac{A_n}{nB_n}= \frac{1}{e}.
$$
By Stirling's formula
$$
 \left(\frac{n}{2}\right)^{\frac{n}{2}}e^{-\frac{n}{2}}\sqrt{\pi n}
\leq
\Gamma\left(1+\frac{n}{2}\right)
\leq \left(\frac{n}{2}\right)^{\frac{n}{2}}e^{-\frac{n}{2}}\sqrt{\pi n}
e^{\frac{1}{6n-12}}
$$
and thus
$$
\sqrt{\frac{n}{2e}} (\pi n)^{\frac{1}{2n}}
\leq
\Gamma\left(1+\frac{n}{2}\right)^{\frac{1}{n}}
\leq\sqrt{\frac{n}{2e}}(\pi n)^{\frac{1}{2n}}e^{\frac{1}{n(6n-12)}}
$$
Moreover,
$$
B_{n}
\geq\frac{  \Gamma(1+\frac{n}{2})^{\frac{2}{n}} } {\Gamma(2) }
\geq\frac{n}{2e}.
$$
(In fact, it can be proved
that $B_n\sim 2e^{c_n^2}\sim  2e^{\frac{n}{2e}}$,  
where we write $a(n)\sim b(n)$, to mean that there are
absolute constants $c_1,c_2>0$ such that $c_1a(n) \leq b(n) \leq c_2\ a(n)$.
But for our purposes the above estimate
is enough.)
Also, for $n\geq 3$
\begin{eqnarray*}
\frac{  \Gamma(1+\frac{n}{2})^{\frac{n-1}{n}} } {\Gamma(1+\frac{n-1}{2}) }
&=&\frac{  \Gamma(1+\frac{n}{2}) } 
{\Gamma(1+\frac{n}{2})^{\frac{1}{n}}\Gamma(1+\frac{n-1}{2}) }   \\
&\leq&\frac{ (\frac{n}{2})^{\frac{n}{2}}e^{-\frac{n}{2}}\sqrt{\pi n}\ e^{\frac{1}{6n-12}}}
{(\pi n)^{\frac{1}{2n}} \sqrt{\frac{n}{2e}} (\frac{n-1}{2})^{\frac{n-1}{2}}e^{-\frac{n-1}{2}}\sqrt{\pi (n-1)}}
\\
&\leq&\frac{ (n)^{\frac{n}{2}}\sqrt{ n}\ } 
{\sqrt{n} (n-1)^{\frac{n-1}{2}}\sqrt{ (n-1)}}  \left( \frac{e}{\pi n}\right)^\frac{1}{2n} \\
&=&  \left( \frac{e}{\pi n}\right)^\frac{1}{2n}  \left(\frac{n}{n-1}\right)^{\frac{n}{2}}
 \leq e^\frac{1}{2n} \left(1+\frac{1}{n-1}\right)^{\frac{n}{2}}\leq\exp\left(\frac{n+1}{2(n-1)}\right)
\leq e.
\end{eqnarray*}
Therefore, by (\ref{centroid2}) and with $c_n =  \big(\Gamma(1+\frac{n}{2})\big)^{\frac{1}{n}}$, 
\begin{eqnarray*}
\lim_{n\to\infty}\frac{A_n}{nB_n}
&=&\lim_{n\to\infty}\frac{2c_n}{n}\left(\frac{1}{\Gamma(\frac{1}{2})B_n}+
c_n\left(1- \frac{\frac{  \Gamma(1+\frac{n}{2})^{\frac{n-1}{n}} } {\Gamma(1+\frac{n-1}{2}) }+
\frac{  \Gamma(1+\frac{n}{2})^{\frac{n}{n}} } {\Gamma(1+\frac{n}{2}) }}{B_n}\right)\right)
\\
&=&\lim_{n\to\infty}\frac{2c_n^2}{n}
=\lim_{n\to\infty} \frac{2}{n}\left(\sqrt{\frac{n}{2e}}\right)^2
=\frac{1}{e}.
\end{eqnarray*}

$\square$
\vskip 4mm

\begin{lemma}\label{prop:gravityM}
Let $a,b \in \mathbb R$, $a>0$, $b>0$ and 
$$
K=\frac{B_2^n}{ \operatorname{vol}_{n} \left(B_2^n \right)^{\frac{1}{n}}}
\hskip 20mm
L=\frac{B_\infty^n}{2}.
$$
Let
$
M_{n}=\operatorname{co}[ (K,-a),(L,b)] 
$
be defined as in (\ref{definition:M}).
Then the center of gravity $g(M_n)$ satisfies
$$
\lim_{n\to \infty} \ g(M_n)=
\frac{1}{e} \ (-a) +\left(1-
\frac{1}{e}\right) \ b.
$$
\end{lemma}
\vskip 3mm
\noindent
{\bf Proof.} 
By symmetry, 
the centroid of $M_n$ is
an element of the
$(n+1)$-axis. Therefore we only need to compute its $(n+1)$-th coordinate.
\par
Instead of $M_n=\mbox{co}[(K,-a),(L,b)]$, we  consider 
$\tilde{M}_n=\mbox{co}[(K,0),(L,a+b)]$ with centroid $\tilde{g}=\tilde{g}_n$. The
centroid
$g_n$ of
$M_n$ is then given by $g_n=\tilde{g}_n -a$.  
\par
We use Lemma \ref{lemma:gravityM} with $c=a+b$ to get
\begin{eqnarray*}
\langle \tilde{g}(\tilde{M}_n),e_{n+1}\rangle
=\frac{a+b}{n+2} \ \frac{\sum_{k=0}^n
(k+1)V_{n-k,k}(K,L)}{\sum_{k=0}^n
V_{n-k,k}(K,L)} .
\end{eqnarray*}
By  Lemma \ref{lemma:mixed volumes}  and the linearity of the
mixed volumes in  each component, we get
\begin{eqnarray}\label{BinomEst4}
\langle \tilde{g}(\tilde{M}_n),e_{n+1}\rangle
&=&\frac{a+b}{n+2}\  \frac{\sum_{k=0}^n
(k+1)\  \mbox{vol}_{n} \left(B_2^n \right)^{\frac{k}{n}} \mbox{vol}_{n-k} \left(B_2^{n-k} \right)}{\sum_{k=0}^n
\mbox{vol}_{n} \left(B_2^n\right)^{\frac{k}{n}} \mbox{vol}_{n-k} \left(B_2^{n-k} \right)} \nonumber\\
&=&\frac{a+b}{n+2}\ \frac{\sum_{k=0}^n
\frac{k+1}{(\Gamma(1+n/2)^{\frac{k}{n}}\  \Gamma(1+(n-k)/2)}}
{\sum_{k=0}^n
\frac{1}{(\Gamma(1+n/2)^{\frac{k}{n}}\  \Gamma(1+(n-k)/2)}}.
\end{eqnarray}
Now we apply Lemma \ref{centroid}.
$\square$
\vskip 3mm

Eventually we will have to investigate expressions of the form
$$
 \mbox{vol}_{n} \left(  \frac{B_2^n } {\mbox{vol}_{n} \left(B_2^n \right)^\frac{1}{n} }\cap  
 \  t\ \frac{B_1^n}{\mbox{vol}_{n} \left(B_1^n \right)^\frac{1}{n} } \right),  \hskip 5mm t \geq 0.
$$
Schechtman, Schmuckenschl\"ager and Zinn established asymptotic
formulas for large $n$ for the volumes of 
$
B_{p}^{n}\cap t B_{q}^{n}
$
\cite{SchSch, SZ, Schmu}.
To do so,  they considered independent random variables
$h_{1}^{p},\dots,h_{n}^{p}$ with density
\begin{equation}\label{PDensity}
\frac{p}{2\Gamma(\frac{1}{p})}e^{-|t|^{p}}
\end{equation}
Here, we need uniform estimates instead of asymptotic ones.
\vskip 3mm

\begin{lemma}
Let $0<p,q<\infty$ and let $h^p$ be a random variable with density
(\ref{PDensity}).
Then
\begin{equation}\label{MomentDens4}
\mathbb E|h^{p}|^{q}
=\frac{p}{2\Gamma(\frac{1}{p})}\int_{-\infty}^{\infty}|t|^{q}e^{-|t|^{p}}dt
=\frac{p}{\Gamma(\frac{1}{p})}\int_{0}^{\infty}t^{q}e^{-t^{p}}dt
=\frac{1}{\Gamma(\frac{1}{p})}\Gamma\left(\frac{q+1}{p}\right)
\end{equation}
\end{lemma}
\par
\noindent
In particular, for $p=1$ and $q=1$ respectively $q=2$ we get
\begin{equation}\label{MomentDens1}
\mathbb E|h^1|^{}=1 \hskip 10mm \mbox{respectively}
\hskip 10mm
\mathbb E|h^1|^{2}=2.
\end{equation}
\par
\noindent
For $p=2$ and  $q=2$ we get
\begin{equation*} 
\mathbb E|h^2|^{2}=\frac{1}{2}.
\end{equation*}

\vskip 3mm

The next lemma treats the Gau{\ss}ian case of a more general statement which we will prove below.

\begin{lemma}\label{ConcentGauss}
Let $(\Omega, \mathbb{P})$ be a probability space.
Let $g_i: \Omega \rightarrow \mathbb{R}$,  $1 \leq i \leq n$,
be independent $N(0,1)$-random
variables. Then, for all
$\gamma>0$
there is $n_{0}$ such that for all $n\geq n_{0}$
\begin{equation}\label{GaussConcentration22}
\mathbb P\left\{\omega\left|\hskip 1mm
\left|
\frac{\frac{1}{n}\sum_{i=1}^{n}|g_{i}(\omega)|}
{\left(\frac{1}{n}\sum_{i=1}^{n}|g_{i}(\omega)|^{2}\right)^{\frac{1}{2}}}
-\sqrt{\frac{2}{\pi}}\right|\leq \gamma
\right.\right\}
\geq\frac{1}{2}
\end{equation}
\end{lemma}
\vskip 3mm

\noindent
{\bf Proof.}
By Chebyshev's inequality.
$\square$
\vskip 3mm

\begin{lemma}\label{ConcentrationRV1}
Let $(\Omega, \mathbb{P})$ be a probability space.
Let  $h_i ^1:  \Omega \rightarrow \mathbb{R}^n$,  $1 \leq i \leq n$, 
be independent
random variables with density $e^{-|t|}$. Then for all 
$\gamma>0$ there is $n_{0}$ such that for all $n\geq n_{0}$
\begin{equation}\label{concentration35}
\mathbb P\left\{\omega\left|\sqrt{2}-\gamma
\leq\frac{(\frac{1}{n}\sum_{i=1}^{n}|h_{i}^1(\omega)|^{2})^{\frac{1}{2}}}
{\frac{1}{n}\sum_{i=1}^{n}|h_{i}^1(\omega)|}
\leq\sqrt{2}+\gamma\right.\right\}
\geq \frac{1}{2}
\end{equation}
\end{lemma}
\vskip 3mm

\noindent
{\bf Proof.}
Again, by Chebyshev's inequality.
$\square$
\vskip 3mm

\begin{lemma}\label{VolGauss1}
Let $(\Omega, \mathbb{P})$ be a probability space.
Let $g_i: \Omega \rightarrow \mathbb{R}$,  $1 \leq i \leq n$,
be independent $N(0,1)$-random
variables.
Then
\begin{equation}\label{VolGauss2}
\mbox{vol}_{n} \left( B_2^n \cap   t\ B_1^{n} \right) 
= n\   \mbox{vol}_{n} \left( B_2^n \right)
\int_0^1 r^{n-1} 
\mathbb P\left(\frac{\sum_{i=1}^n |g_i|}{(\sum_{i=1}^n
|g_i|^2)^\frac{1}{2}}
 \leq \frac{t}{r} \right)\ dr
\end{equation}
\end{lemma}
\vskip 4mm 
\noindent
Again, we state this lemma and its proof in the Gau{\ss}ian case to illustrate the ideas.   Below, we will prove a more general statement. 
A more general situation has also been explored in
\cite{SZ}. 
\par
\noindent
Clearly, formula (\ref{VolGauss2}) is equivalent to
\begin{equation}\label{VolGauss3}
\mbox{vol}_{n} \left( B_2^n \cap  t\ B_1^{n} \right) 
= n \ \mbox{vol}_{n} \left(  B_2^n \right)\int_0^1 r^{n-1} \ 
\mathbb P\bigg(
\frac{\frac{1}{n}\sum_{i=1}^{n}|g_i|}
{(\frac{1}{n}\sum_{i=1}^n|g_i|^2)^\frac{1}{2}}
 \leq \frac{t}{r\sqrt{n}} \bigg)\ dr
\end{equation}
and to
\begin{eqnarray}\label{VolGauss4}
&& \mbox{vol}_{n} \left( \frac{B_2^n}{\mbox{vol}_{n} \left(B_{2}^{n}\right)^{1/n}} \cap  s\
\frac{B_1^{n}}{\mbox{vol}_{n} \left(B_{1}^{n}\right)^{1/n}}
\right)   \nonumber \\
&&= n  \int_0^1 r^{n-1} \ 
\mathbb P\bigg(
\frac{\frac{1}{n}\sum_{i=1}^{n}|g_i|}
{(\frac{1}{n}\sum_{i=1}^n|g_i|^2)^\frac{1}{2}}
 \leq \frac{s}{r\sqrt{n}}
\left(\frac{\mbox{vol}_{n} \left( B_2^n\right)}{\mbox{vol}_{n} \left( B_1^n\right)}\right)^{\frac{1}{n}} \bigg)\
dr.
\end{eqnarray}
(\ref {VolGauss4}) follows from (\ref {VolGauss3}) with the substitution $s=t \left(\frac{\mbox{vol}_{n} \left(B_{1}^{n}\right)}{\mbox{vol}_{n} \left(B_{2}^{n}\right)}\right)^{1/n}$.

\vskip 3mm
\noindent
{\bf Proof of Lemma \ref{VolGauss1}.}
Also in the case of a convex body that is not centrally symmetric we use
$$
\|\xi\|_{K}
=\inf\left\{\left.\frac{1}{\rho}\right|\rho\xi\in K\right\}\ .
$$
Using polar coordinates we get for any convex body  $K \subseteq B^n_2$ 
\begin{eqnarray*}
\mbox{vol}_{n} \left( K \right) 
&=&\frac{1}{n}\int_{\partial B_{2}^{n}}\|\xi\|^{-n}_K d\sigma_{n-1}(\xi)
\\
&=&\int_{\partial B_{2}^{n}}
\int_{0}^{\|\xi\|_{K}^{-1}}r^{n-1}drd\sigma_{n-1}(\xi) \\
&=&\int_0^{1} r^{n-1}\operatorname{vol}_{n-1}\left(
\partial B_{2}^{n}\cap \frac{1}{r}K\right)\ dr.
\end{eqnarray*}
Hence for $K=B_{2}^{n}\cap t B_{1}^{n}$
\begin{eqnarray*}
\mbox{vol}_{n} \left( B_{2}^{n}\cap t B_{1}^{n} \right)
&=&\int_0^{1} r^{n-1}\operatorname{vol}_{n-1}\left(
\partial B_{2}^{n}\cap \frac{1}{r}(B_{2}^{n}\cap t B_{1}^{n})\right)\
dr   \\
&=&\int_0^{1} r^{n-1}\operatorname{vol}_{n-1}\left(
\partial B_{2}^{n}\cap \frac{t}{r} B_{1}^{n}\right)\ dr
\end{eqnarray*}
Let $\gamma_{n}$ be the normalized Gauss measure on $\mathbb R^{n}$.
Then
$$
\operatorname{vol}_{n-1}(\partial B_{2}^{n})\cdot
\gamma_{n}\left(\frac{\sum_{i=1}^n |x_i|}{(\sum_{i=1}^n
|x_i|^2)^\frac{1}{2}}
 \leq \frac{t}{r} \right)
=\operatorname{vol}_{n-1}\left(
\partial B_{2}^{n}\cap \frac{t}{r} B_{1}^{n}\right)
$$
Indeed,
\begin{eqnarray*}
\left\{x\left|\frac{\sum_{i=1}^n |x_i|}{\left(\sum_{i=1}^n
|x_i|^2\right)^\frac{1}{2}}
 \leq \frac{t}{r} \right.\right\}
&=&\left\{sx\left|s>0, \|x\|_{2}=1,
\|x\|_{1}\leq\frac{t}{r}\right.\right\}   \\
&=&\left\{sx\left|s>0,x\in\partial B_{2}^{n}\cap \frac{t}{r}B_{1}^{n}
\right.\right\}
\end{eqnarray*}
Therefore, using polar coordinates
\begin{eqnarray*}
\gamma_n \left(\left\{x\left|\frac{\sum_{i=1}^n |x_i|}{\left(\sum_{i=1}^n
|x_i|^2\right)^\frac{1}{2}}
 \leq \frac{t}{r} \right.\right\} \right)
&=& \gamma_n \left( \left\{sx\left|s>0,x\in\partial B_{2}^{n}\cap \frac{t}{r}B_{1}^{n}
\right.\right\} \right) \\
&= &\frac{1}{(2 \pi)^\frac{n}{2}} \int_{s >0} \int_{\partial B_{2}^{n}\cap \frac{t}{r}B_{1}^{n}} s^{n-1} e^{-\frac{s^2}{2}} dx \  ds\\
&=& \frac{\Gamma\left(\frac{n}{2}\right)}{2 \pi ^\frac{n}{2}}\mbox{vol}_{n-1} \left(  \partial B_{2}^{n}\cap \frac{t}{r}B_{1}^{n} \right).
\end{eqnarray*}
Thus
$$
\mbox{vol}_{n} \left( B_{2}^{n}\cap t B_{1}^{n} \right)
=\operatorname{vol}_{n-1}(\partial B_{2}^{n})\int_0^{1}
r^{n-1}\gamma_{n}\left(\frac{\sum_{i=1}^n |x_i|}{(\sum_{i=1}^n
|x_i|^2)^\frac{1}{2}}
 \leq \frac{t}{r} \right)\ dr.
$$
Since the coordinate functionals are independent $N(0,1)$-random
variables with respect to the measure $\gamma_{n}$ we have
established (\ref{VolGauss2}).
$\square$
\vskip 3mm

\vskip 3mm\noindent
\begin{lemma} \label{}
Let ${\mathbb P}_p$ be the probability on $\mathbb{R}^n$ with density $f_p:\mathbb{R}^n \to 
\mathbb{R}$
defined by
$$f_p(x)=f_p(x_1,\dots,x_n) =
\frac{1}{\big(2\Gamma(1+\frac{1}{p})\big)^n} e^{-\sum_{i=1}^n |x_i|^p}.$$
Then for every starshaped body $K$, one has 
$$
\operatorname{vol}_{n}\left(B_p^n \cap K\right)
= n  \  \operatorname{vol}_{n}\left(B_p^n \right)
 \int_0^1 r^{n-1} \mathbb{P}_p\left(\left\{x\in\mathbb{R}^n\left|
\frac{\|x\|_K} {\|x\|_p} \leq \frac{1}{r}\right.\right\}\right)dr. 
$$
\end{lemma}

\vskip 3mm
\noindent
{\bf Proof.}
As $\mbox{vol}_{n}\left(B_p^n\right) =\big(2\Gamma(1+\frac{1}{p})\big)^n
\big(\Gamma(1+\frac{n}{p})\big)^{-1}$, it follows that
 \begin{eqnarray*}
&&n \  \mbox{vol}_{n}\left(B_p^n\right) \int_0^1 r^{n-1} {\mathbb P}_p
\big(\big\{x\in \mathbb {R}^n;
\frac{\|x\|_K} {\|x\|_p} \le \frac{1}{r}\big\}\big)dr   \\
&&=
\frac{n}{\Gamma(1+\frac{n}{p})}
\int_0^1\Big(\int_{ \{x: r\|x\|_K\le\|x\|_p \}} 
e^{-\|x\|_p^p}dx\Big) r^{n-1}dr
\\
&&=\frac{n}{\Gamma(1+\frac{n}{p})}\int_{\mathbb{R}^n}
\Big( \int_0^{\min(1,\frac{\|x\|_p}{\|x\|_K})}r^{n-1}dr\Big) 
e^{-\|x\|_p^p}dx
\\
&&=\frac{1}{\Gamma(1+\frac{n}{p})}
\int_{\mathbb{R}^n}\frac{1}{\big(\max(1,\frac{\|x\|_K}{\|x\|_p})\big)^n}
e^{-\|x\|_p^p}dx
\end{eqnarray*}
We pass to polar coordinates. $\sigma_{n-1}$ denotes the (non-normalized)
surface measure on $S^{n-1}$.
\begin{eqnarray*}
&&n \  \mbox{vol}_{n}\left(B_p^n\right) \int_0^1 r^{n-1} {\mathbb P}_p
\big(\big\{x\in \mathbb {R}^n;
\frac{\|x\|_K} {\|x\|_p} \le \frac{1}{r}\big\}\big)dr   \\
&& =\frac{1}{\Gamma(1+\frac{n}{p})}\int_{\theta\in S_{n-1}}
\Big(\int_0^{+\infty} e^{-r^p\|\theta\|_p^{p}} r^{n-1} dr\Big)
\frac{1}{\big(\max(1,\frac{\|\theta\|_K}{\|\theta\|_p})\big)^n}d\sigma_{n-1}(\theta)
\\
&&= \frac{1}{\Gamma(1+\frac{n}{p})}\frac{1}{p} 
\int_0^{+\infty} e^{-s} s^{\frac{n}{p} -1} ds
\int_{\theta\in S_{n-1}}\frac{1}{\|\theta\|_p^n}
\frac{1}{\big(\max(1,\frac{\|\theta\|_K}{\|\theta\|_p})\big)^n}d\sigma_{n-1}(\theta)
\\
&&= \frac{1}{n}\int_{\theta\in
S_{n-1}}\frac{1}{\big(\max(\|\theta\|_p, \|\theta\|_K)\big)^n}
d\sigma_{n-1}(\theta)   \\
&&= \operatorname{vol}_{n}\left(K\cap B_p^n \right)
\end{eqnarray*}
$\square$

\vskip 3mm
\noindent
\begin{cor} \label{corvol}
Let $1 \leq p,q, < \infty$ and let $s \geq 0$. 
Let  $h_i ^p:  \Omega \rightarrow \mathbb{R}^n$,  $1 \leq i \leq n$, 
be independent
random variables with density 
$\frac{p}{\Gamma(1/p)}e^{-|t|^{p}}$ and $h=(h_{1}^p,\dots,h_{n}^p)$.
Then
\begin{eqnarray*}
&&\operatorname{vol}_{n}\left(B_p^n\cap s
B_q^n\right)
\\
&&=n \ \operatorname{vol}_{n}\left(B_p^n \right) \int_0^1 r^{n-1} 
\mathbb{P}_p\left(\left\{ (x_1,\dots,x_n)\in \mathbb{R}^n\left|
\frac{\big(\frac{1}{n} \sum_{i=1}^n |x_i|^q\big)^{\frac{1}{q}} }
     {\big(\frac{1}{n} \sum_{i=1}^n |x_i|^p\big)^{\frac{1}{p}} }
\le  n^{\frac{1}{p}-\frac{1}{q}} 
 \frac{s}{r}\right.\right\}
\right) dr    \\
&&=
n \int_{0}^{1}r^{n-1}
\mathbb P
\left\{
\omega\left|\frac{\|h^{p}(\omega)\|_{q}}{\|h^{p}(\omega)\|_{p}}\leq \frac{s}{r} 
\right.\right\} dr
\end{eqnarray*}
\end{cor}
\vskip 3mm

\begin{lemma}\label{VolIntersect51}
For all $\gamma>0$ there is $n_{0}$ such that for all $n\geq n_{0}$
and all $t$ with
\newline
$t\geq\left(\sqrt{\frac{2}{\pi}}+\gamma\right)  \sqrt{n} \  
\frac{\operatorname{vol}_{n}\left(B_1^n\right)^{\frac{1}{n}}}
{\operatorname{vol}_{n}\left(B_2^n\right)^{\frac{1}{n}}}$
we have
$$
\frac{1}{2}
\leq
\operatorname{vol}_{n}\left(\frac{B_2^n}{\operatorname{vol}_{n}\left(B_{2}^{n}\right)^{1/n}} \cap  t\
\frac{B_1^{n}}{\operatorname{vol}_{n}\left(B_{1}^{n}\right)^{1/n}}
\right)  
\leq 1
\ .
$$
\end{lemma}
\vskip 3mm

\noindent
{\bf Proof.}
By (\ref{VolGauss3})
$$
\mbox{vol}_{n}\left( B_2^n \cap  t\ B_1^{n} \right)
= n \  \mbox{vol}_{n}\left( B_2^n \right) \int_0^1 r^{n-1} \ 
\mathbb P\bigg(
\frac{\frac{1}{n}\sum_{i=1}^{n}|g_i|}
{(\frac{1}{n}\sum_{i=1}^n|g_i|^2)^\frac{1}{2}}
 \leq \frac{t}{r\sqrt{n}} \bigg)\ dr
 \ .
$$
This gives
\begin{equation}\label{volfor}
\mbox{vol}_{n}\left( \frac{B_2^n}{\mbox{vol}_{n}\left(B_2^n\right)^{\frac{1}{n}}} 
\cap  t\frac{\ B_1^{n}}{\mbox{vol}_{n}\left( B_2^n\right)^{\frac{1}{n}}} \right)
= n \int_0^1 r^{n-1} \ 
\mathbb P\bigg(
\frac{\frac{1}{n}\sum_{i=1}^{n}|g_i|}
{(\frac{1}{n}\sum_{i=1}^n|g_i|^2)^\frac{1}{2}}
 \leq \frac{t}{r\sqrt{n}} \bigg)\ dr.
\end{equation}
We substitute
$$
t=s\frac{\mbox{vol}_{n}\left(B_2^n\right)^{\frac{1}{n}}}{\mbox{vol}_{n}\left(B_1^n\right)^{\frac{1}{n}}}
$$
and obtain
\begin{eqnarray*}
&& \mbox{vol}_{n}\left(\frac{B_2^n}{\mbox{vol}_{n}\left(B_2^n\right)^{\frac{1}{n}}} 
\cap  s\frac{\ B_1^{n}}{\mbox{vol}_{n}\left(B_1^n\right)^{\frac{1}{n}}} \right) \\
&& = n \int_0^1 r^{n-1} \ 
\mathbb P\left(
\frac{\frac{1}{n}\sum_{i=1}^{n}|g_i|}
{(\frac{1}{n}\sum_{i=1}^n|g_i|^2)^\frac{1}{2}}
 \leq \frac{s}{r\sqrt{n}}\frac{\mbox{vol}_{n}\left(B_2^n\right)^{\frac{1}{n}}}{\mbox{vol}_{n}\left(B_1^n\right)^{\frac{1}{n}}} \right)\ dr.
\end{eqnarray*}
By Lemma \ref{ConcentGauss}, for every $\gamma>0$ there is $n_{0}$ such that for all $n\geq n_{0}$
$$
\mathbb P\left\{\omega\left|\hskip 1mm
\left|
\frac{\frac{1}{n}\sum_{i=1}^{n}|g_{i}(\omega)|}
{\left(\frac{1}{n}\sum_{i=1}^{n}|g_{i}(\omega)|^{2}\right)^{\frac{1}{2}}}
-\sqrt{\frac{2}{\pi}}\right|\leq \gamma
\right.\right\}
\geq\frac{1}{2}
$$
Therefore, for every $\gamma>0$ there is $n_{0}$ such that for all $n\geq n_{0}$
and all $s$ with
$$
\sqrt{\frac{2}{\pi}}+\gamma
\leq\frac{s}{\sqrt{n}}  \   \frac{\mbox{vol}_{n}\left(B_2^n\right)^{\frac{1}{n}}}{\mbox{vol}_{n}\left(B_1^n\right)^{\frac{1}{n}}}
$$
or, equivalently,
$$
\left(\sqrt{\frac{2}{\pi}}+\gamma\right)\sqrt{n} \  
\frac{\mbox{vol}_{n}\left(B_1^n\right)^{\frac{1}{n}}}{\mbox{vol}_{n}\left(B_2^n\right)^{\frac{1}{n}}}
\leq s
$$
we have
$$
\mbox{vol}_{n}\left( \frac{B_2^n}{\mbox{vol}_{n}\left(B_2^n\right)^{\frac{1}{n}}} 
\cap  s\frac{\ B_1^{n}}{\mbox{vol}_{n}\left(B_1^n\right)^{\frac{1}{n}}} \right)
\geq\frac{1}{2}.
$$
The other inequality is obvious from (\ref{volfor}).
$\square$
\vskip 3mm

\begin{lemma}\label{VolIntersect22}
For all $\gamma>0$ there is $n_{0}$ such that for all $n\geq n_{0}$
and all $s$ with
$$
s \geq  \frac{\sqrt{2}+\gamma}  
{\sqrt{n}} \ 
\frac{\operatorname{vol}_{n}\left(B_{2}^{n}\right)^{\frac{1}{n}}}
{\operatorname{vol}_{n}\left(B_{1}^{n}\right)^{\frac{1}{n}}}
\sim\sqrt{\frac{\pi}{e}}
$$
we have
$$
\frac{1}{2}
\leq
\operatorname{vol}_{n}\left(\frac{B_{1}^{n}}
{\operatorname{vol}_{n}\left(B_{1}^{n}\right)^{\frac{1}{n}}}\cap
s\frac{B_{2}^{n}}{\operatorname{vol}_{n}\left(B_{2}^{n}\right)^{\frac{1}{n}}}\right)\leq 1
\ .
$$
\end{lemma}
\vskip 3mm

\noindent
{\bf Proof.} By Corollary \ref{corvol}, for $q=2$ and $p=1$
$$
\frac{\operatorname{vol}_{n}(B_{1}^{n}\cap tB_{2}^{n})}
{\operatorname{vol}_{n}(B_{1}^{n})}
=
n \int_{0}^{1}r^{n-1}
\mathbb P
\left\{
\omega\left|\frac{\|h^{1}(\omega)\|_{2}}{\|h^{1}(\omega)\|_{1}}\leq \frac{t}{r}
\right.\right\}
dr\ .
$$
We put
$$
t=s\frac{\operatorname{vol}_{n}\left(B_{1}^{n}\right)^{\frac{1}{n}}}
{\operatorname{vol}_{n}\left(B_{2}^{n}\right)^{\frac{1}{n}}}
$$
and obtain
\begin{eqnarray}\label{VolIntersec11}
 \operatorname{vol}_{n}\left(\frac{B_{1}^{n}}
{\operatorname{vol}_{n}\left(B_{1}^{n}\right)^{\frac{1}{n}}}\cap
s\frac{B_{2}^{n}}{\mbox{vol}_{n}\left(B_{2}^{n}\right)^{\frac{1}{n}}}\right)
= 
 n \int_{0}^{1}r^{n-1}\mathbb P\left(
\frac{\|h^{1}(\omega)\|_{2}}{\|h^{1}(\omega)\|_{1}}\leq\frac{s}{r}
\frac{\operatorname{vol}_{n}\left(B_{1}^{n}\right)^{\frac{1}{n}}}
{\operatorname{vol}_{n}\left(B_{2}^{n}\right)^{\frac{1}{n}}}
\right)dr.
\end{eqnarray}
By Lemma \ref{ConcentrationRV1},  for every $\gamma>0$ there is $n_{0}$ such that for all
$n\geq n_{0}$
$$
\mathbb P\left\{\omega\left|
\frac{(\frac{1}{n}\sum_{i=1}^{n}|h_{i}^{1}(\omega)|^{2})^{\frac{1}{2}}}
{\frac{1}{n}\sum_{i=1}^{n}|h_{i}^{1}(\omega)|}
\leq\sqrt{2}+\gamma\right.\right\}
\geq \frac{1}{2}
$$
Therefore, for all $r$ with $0<r\leq 1$ and all $s$ with
$$
s \  \  \frac{\mbox{vol}_{n}\left(B_{1}^{n}\right)^{\frac{1}{n}}}{\mbox{vol}_{n}\left(B_{2}^{n}\right)^{\frac{1}{n}}}
\geq \frac{\sqrt{2}+\gamma}  
{\sqrt{n}}
$$
or, equivalently,
$$
s \geq \frac{\sqrt{2}+\gamma}  
{\sqrt{n}}\  
\frac{\mbox{vol}_{n}\left(B_{2}^{n}\right)^{\frac{1}{n}}}{\mbox{vol}_{n}\left(B_{1}^{n}\right)^{\frac{1}{n}}}
$$
we have 
$$
\mathbb P\left\{\omega\left|
\frac{(\frac{1}{n}\sum_{i=1}^{n}|h_{i}^{1}(\omega)|^{2})^{\frac{1}{2}}}
{\frac{1}{n}\sum_{i=1}^{n}|h_{i}^{1}(\omega)|}
\leq s\frac{\mbox{vol}_{n}\left(B_{1}^{n}\right)^{\frac{1}{n}}}{\mbox{vol}_{n}\left(B_{2}^{n}\right)^{\frac{1}{n}}}\right.\right\}
\geq \frac{1}{2}.
$$
By (\ref{VolIntersec11})
$$
\frac{1}{2}
\leq 
\operatorname{vol}_{n}\left(\frac{B_{1}^{n}}{\mbox{vol}_{n}\left(B_{1}^{n}\right)^{\frac{1}{n}}}\cap
s\frac{B_{2}^{n}}{\mbox{vol}_{n}\left(B_{2}^{n}\right)^{\frac{1}{n}}}\right)
$$
provided that
$$
s\geq \frac{\sqrt{2}+\gamma}  
{\sqrt{n}}\  
\frac{\mbox{vol}_{n}\left(B_{2}^{n}\right)^{\frac{1}{n}}}{\mbox{vol}_{n}\left(B_{1}^{n}\right)^{\frac{1}{n}}}.
$$
Again, the other inequality follows by (\ref{VolIntersec11}).
$\square$
\vskip 3mm

\noindent
Let $a>0, b>0$ be real numbers.
By (\ref{Mo(s)}), for $K= \frac{B^n_2}{\mbox{vol}_{n}\left(B_{2}^{n}\right)^{\frac{1}{n}}} $, 
$L= \frac{B^n_\infty}{2}$, for all $s$ with $-\frac{1}{a}\leq s\leq \frac{1}{b}$
\begin{eqnarray}\label{FormulaM11}
M_n^\circ(s)=(1+s \ a) \mbox{vol}_{n}\left( B_2^n\right)^{\frac{1}{n}}B_{2}^{n} 
\cap (1-s \ b) 2 B_{1}^{n}
\end{eqnarray}
which is the same as
\begin{eqnarray}\label{FormulaM1}
&M_n^\circ(s)
= \nonumber \\
&(1+s \ a) \mbox{vol}_{n}\left( B_2^n \right)^{\frac{2}{n}}
\left(\frac{B_{2}^{n}}{\mbox{vol}_{n}\left(B_2^n\right)^{\frac{1}{n}}} \ 
\cap\  2 \  \frac{(1-s \ b)}{(1+s \ a)}
\frac{\mbox{vol}_{n}\left(B_1^n \right)^{\frac{1}{n}}}{\mbox{vol}_{n}\left(B_2^n\right)^{\frac{2}{n}}}\frac{B_{1}^{n}}
{\mbox{vol}_{n}\left(B_1^n\right)^{\frac{1}{n}}} \right).
\end{eqnarray}
\vskip 3mm

\begin{lemma}\label{Kegel1}
Let $a>0, b>0$ be real numbers.
Then, for all $\gamma>0$ there is
$n_{0}$ such that for all $n\geq n_{0}$ and all $s$ with
$$
s \geq 
- \frac{ \ \sqrt{\pi e}- 2}
{a \  \ \sqrt{\pi e} + 2 b } +\gamma 
$$
we have
\begin{eqnarray*}
2^{n-1} \  (1-s \ b)^{n}\operatorname{vol}_{n}(B_{1}^{n})
\leq
\mbox{vol}_{n}\left( M_n^{\circ}(s)
\right)  
\leq2^{n}(1-s \ b)^{n}\operatorname{vol}_{n}(B_{1}^{n}).
\end{eqnarray*}
\end{lemma}
\vskip 3mm

\noindent
{\bf Proof.}
We have
$$
M_n^\circ(s)=2(1-s \ b)\left(
\frac{(1+s \ a)  \mbox{vol}_{n}\left(B_2^n \right) ^{\frac{1}{n}}}{2(1-s \ b)}B_{2}^{n} 
\cap B_{1}^{n}\right)
$$
Therefore
\begin{eqnarray*}
&&\operatorname{vol}_{n}(M_n^\circ(s)) 
=2^{n}(1-s \ b)^{n}
\operatorname{vol}_{n}\left(
\frac{(1+s \ a)  \mbox{vol}_{n}\left( B_2^n \right)^{\frac{1}{n}}}{2(1-s \ b)}B_{2}^{n} 
\cap B_{1}^{n}\right)   \\
&&=2^{n}(1-s \ b)^{n}\operatorname{vol}_{n}(B_{1}^{n})
\operatorname{vol}_{n}\left(
\frac{(1+s \ a)}{2(1-s \ b)}
\frac{ \mbox{vol}_{n}\left( B_2^n \right)^{\frac{2}{n}}}{ \mbox{vol}_{n}\left(B_{1}^{n}\right)^{\frac{1}{n}}}
\frac{B_{2}^{n}}{ \mbox{vol}_{n}\left( B_2^n\right)^{\frac{1}{n}}}
\cap \frac{B_{1}^{n}}{ \mbox{vol}_{n}\left(B_{1}^{n}\right)^{\frac{1}{n}}}\right).
\end{eqnarray*}
We get immediately
$$
\operatorname{vol}_{n}(M_n^\circ(s)) 
\leq2^{n}(1-s \ b)^{n}\operatorname{vol}_{n}(B_{1}^{n}).
$$
We show now the opposite inequality.
By Lemma \ref{VolIntersect22}, for all $\gamma>0$ there is $n_{0}$ such that for all  $n\geq n_{0}$
\begin{equation}\label{absch1}
\frac{1}{2}
\leq
\operatorname{vol}_{n}\left(\frac{B_{1}^{n}}{\mbox{vol}_{n}\left(B_{1}^{n}\right)^{\frac{1}{n}}}\cap
\frac{(1+s \ a)}{2(1-s \ b)}
\frac{\mbox{vol}_{n}\left( B_2^n\right)^{\frac{2}{n}}}{\mbox{vol}_{n}\left(B_{1}^{n}\right)^{\frac{1}{n}}}
\frac{B_{2}^{n}}{\mbox{vol}_{n}\left(B_{2}^{n}\right)^{\frac{1}{n}}}\right)\leq 1,
\end{equation}
provided that
$$
\frac{(1+s \ a)}{2(1-s \ b)} \  
\frac{\mbox{vol}_{n}\left( B_2^n\right)^{\frac{2}{n}}}{\mbox{vol}_{n}\left(B_{1}^{n}\right)^{\frac{1}{n}}}
\geq \frac{\sqrt{2}+\gamma}{\sqrt{n}} \  
\frac{\mbox{vol}_{n}\left(B_{2}^{n}\right)^{\frac{1}{n}}}{\mbox{vol}_{n}\left(B_{1}^{n}\right)^{\frac{1}{n}}}
$$
or
$$
\frac{1+s \ a}{1-s \ b}
\geq 2 \ 
\frac{(\sqrt{2}+\gamma)}{\sqrt{n} \  \mbox{vol}_{n}\left(B_{2}^{n} \right)^{\frac{1}{n}}}.
$$
This means
$$
s \geq - \frac{\sqrt{n} \  \mbox{vol}_{n}\left(B_{2}^{n}\right)^{\frac{1}{n}}- 2 (\sqrt{2}+\gamma)}{a \ \sqrt{n} \  \mbox{vol}_{n}\left(B_{2}^{n}\right)^{\frac{1}{n}} + 2 b (\sqrt{2}+\gamma)} .
$$
Since
$$
\operatorname{vol}_{n}(B_{2}^{n})
=\frac{\pi^{\frac{n}{2}}}{\Gamma(\frac{n}{2}+1)}
$$
and
$$
k^{k}e^{-k}\sqrt{2\pi k}
\leq
\Gamma(k+1)
\leq k^{k}e^{-k}\sqrt{2\pi k}e^{\frac{1}{6(k-2)}}
$$
we have
\begin{equation}\label{Kegel20}
\frac{\sqrt{2\pi e}}{\sqrt{n}(\pi n)^{\frac{1}{2n}}e^{\frac{1}{6n(n-2)}}}
\leq
\operatorname{vol}_{n}(B_{2}^{n})^{\frac{1}{n}}
\leq\frac{\sqrt{2\pi e}}{\sqrt{n}(\pi n)^{\frac{1}{2n}}}.
\end{equation}
Thus (\ref{absch1}) holds, provided that 
$$
s \geq - \frac{ \  \frac{\sqrt{2\pi e}}{(\pi n)^{\frac{1}{2n}}}- 2 (\sqrt{2}+\gamma)}
{a \  \ \frac{\sqrt{2\pi e}}{(\pi n)^{\frac{1}{2n}}e^{\frac{1}{6n(n-2)}}} + 2 b (\sqrt{2}+\gamma)}
=- \frac{ \ \sqrt{2\pi e}- 2 (\sqrt{2}+\gamma)(\pi n)^{\frac{1}{2n}}}
{a \  \ \frac{\sqrt{2\pi e}}{e^{\frac{1}{6n(n-2)}}} + 2 b (\sqrt{2}+\gamma)(\pi n)^{\frac{1}{2n}}}  .
$$
Now we pass to a new $\gamma$ and obtain: For every $\gamma>0$ there is a
$n_{0}$ such that for all $n\geq n_{0}$ and all $s$ with
$$
s \geq 
- \frac{ \ \sqrt{2\pi e}- 2 \sqrt{2}}
{a \  \ \sqrt{2\pi e} + 2 b \sqrt{2}} +\gamma
=- \frac{ \ \sqrt{\pi e}- 2}
{a \  \ \sqrt{\pi e} + 2 b } +\gamma.
$$
we have
\begin{eqnarray} \label{18ii}
\frac{1}{2}
\leq
\operatorname{vol}_{n}\left(\frac{B_{1}^{n}}{\mbox{vol}_{n}\left(B_{1}^{n}\right)^{\frac{1}{n}}}\cap
\frac{(1+s \ a)}{2 \ (1-s \ b)}
\frac{\mbox{vol}_{n}\left( B_2^n \right) ^{\frac{2}{n}}}{\mbox{vol}_{n}\left(B_{1}^{n} \right)^{\frac{1}{n}}}
\frac{B_{2}^{n}}{\mbox{vol}_{n}\left(B_{2}^{n} \right)^{\frac{1}{n}}}\right)\leq 1.
\end{eqnarray}
$\square$
\vskip 3mm

\noindent
\begin{lemma}\label{Kegel2}
Let $a>0$ and $b>0$ real numbers. Then, for all $\gamma>0$ there is $n_{0}$
such that for all $n\geq n_{0}$ and all $s$ with
$$
s\leq-\frac{\sqrt{e}-1}{b
+a\sqrt{e}}
-\gamma
$$
we have
$$
 \frac{1}{2}(1+s \ a)^{n}\operatorname{vol}_{n}(B_{2}^{n})^{2}
\leq
\mbox{vol}_{n}\left( M_n^{\circ}(s)
\right)  
\leq(1+s \ a)^{n}\operatorname{vol}_{n}(B_{2}^{n})^{2}
$$
\end{lemma}
\vskip 3mm

\noindent
{\bf Proof.}
By (\ref{FormulaM1}), for all $s$ with $-\frac{1}{a}\leq s\leq \frac{1}{b}$
\begin{equation*} 
M_n^\circ(s)
=(1+s \ a) \mbox{vol}_{n}\left( B_2^n \right)^{\frac{2}{n}}
\left(\frac{B_{2}^{n}}{\mbox{vol}_{n}\left(B_2^n\right)^{\frac{1}{n}}} \ 
\cap\  2 \  \frac{(1-s \ b)}{(1+s \ a)}
\frac{\mbox{vol}_{n}\left(B_1^n \right)^{\frac{1}{n}}}{\mbox{vol}_{n}\left(B_2^n\right)^{\frac{2}{n}}}\frac{B_{1}^{n}}
{\mbox{vol}_{n}\left(B_1^n\right)^{\frac{1}{n}}} \right).
\end{equation*}
The right hand side inequality follows immediately. We show now the left hand side
inequality.
By Lemma \ref{VolIntersect51}, 
for all $\gamma>0$ there is $n_{0}$ such that for all $n\geq n_{0}$
and all $t$ with
\newline
$t\geq\left(\sqrt{\frac{2}{\pi}}+\gamma\right)  \sqrt{n} \  
\frac{\mbox{vol}_{n}\left(B_1^n\right)^{\frac{1}{n}}}{\mbox{vol}_{n}\left(B_2^n\right)^{\frac{1}{n}}}$
we have
$$
\frac{1}{2}
\leq
\mbox{vol}_{n}\left(\frac{B_2^n}{\mbox{vol}_{n}\left(B_{2}^{n}\right)^{1/n}} \cap  t\
\frac{B_1^{n}}{\mbox{vol}_{n}\left(B_{1}^{n}\right)^{1/n}}
\right)  
\leq 1.
$$
Therefore, 
$$
\frac{1}{2}
\leq
\mbox{vol}_{n}\left(\frac{B_2^n}{\mbox{vol}_{n}\left(B_{2}^{n}\right)^{1/n}} \cap  
2 \  \left(\frac{1-s \ b}{1+s \ a}\right)
\frac{\mbox{vol}_{n}\left(B_1^n \right)^{\frac{1}{n}}}{\mbox{vol}_{n}\left(B_2^n\right)^{\frac{2}{n}}}
\frac{B_1^{n}}{\mbox{vol}_{n}\left(B_{1}^{n}\right)^{1/n}}
\right)  
\leq 1
$$
provided that
$$
\left(\sqrt{\frac{2}{\pi}}+\gamma\right)  \sqrt{n} \  
\frac{\mbox{vol}_{n}\left(B_1^n\right)^{\frac{1}{n}}}{\mbox{vol}_{n}\left(B_2^n\right)^{\frac{1}{n}}}
\leq
2 \  \left(\frac{1-s \ b}{1+s \ a}\right)
\frac{\mbox{vol}_{n}\left(B_1^n \right)^{\frac{1}{n}}}{\mbox{vol}_{n}\left(B_2^n\right)^{\frac{2}{n}}}.
$$
The latter inequality is equivalent to
$$
\left(\sqrt{\frac{2}{\pi}}+\gamma\right)  \sqrt{n} \  
\mbox{vol}_{n}\left(B_2^n\right)^{\frac{1}{n}}
\leq
2 \  \frac{1-s \ b}{1+s \ a}.
$$
By (\ref{Kegel20}),
the above inequality holds if
$$
\left(\sqrt{\frac{2}{\pi}}+\gamma\right) 
\sqrt{\frac{\pi e}{2}}
\leq
 \frac{1-s \ b}{1+s \ a}(\pi n)^{\frac{1}{2n}}.
$$
This is equivalent to
$$
s\leq-\frac{\left(\sqrt{\frac{2}{\pi}}+\gamma\right) 
\sqrt{\frac{\pi e}{2}}-(\pi n)^{\frac{1}{2n}}}{b(\pi n)^{\frac{1}{2n}}
+a\left(\sqrt{\frac{2}{\pi}}+\gamma\right) 
\sqrt{\frac{\pi e}{2}}}\ .
$$
We pass to a new $\gamma$
$$
s\leq-\frac{\sqrt{e}-1}{b
+a\sqrt{e}}
-\gamma
$$
$\square$
\vskip 3mm

\begin{lemma}\label{ConcentrationM}
Let 
$$
s_{0}=\frac{1-\sqrt{e}}{b+a\sqrt{e}}
\hskip 30mm
s_{1}
=
\frac{2-\sqrt{\pi e}}{a\sqrt{\pi e}+2b}
$$
\newline
(i) For every $\gamma>0$ there is $n_{0}$ such that for all
$n$ with $n\geq n_{0}$
$$
\int_{-\frac{1}{a}}^{\frac{1}{b}}  \operatorname{vol}_{n}(M_{n}^\circ(s) \ ds
\leq(1+\gamma)\int_{s_{0}-\gamma}^{s_{1}+\gamma}  
\operatorname{vol}_{n}(M_{n}^\circ(s)) \ ds
$$
(ii) For every $\gamma>0$ there is $n_{0}$ such that for all
$n$ with $n\geq n_{0}$
$$
\int_{-\frac{1}{a}}^{s_{0}-\gamma} |s |
\operatorname{vol}_{n}(M_{n}^\circ(s) )\ ds
\leq\gamma
\int_{-\frac{1}{a}}^{\frac{1}{b}}  \operatorname{vol}_{n}(M_{n}^\circ(s)) \ ds
$$
and
$$
\int_{s_{1}+\gamma}^{\frac{1}{b}} |s |
\operatorname{vol}_{n}(M_{n}^\circ(s) )\ ds
\leq\gamma
\int_{-\frac{1}{a}}^{\frac{1}{b}}  \operatorname{vol}_{n}(M_{n}^\circ(s)) \ ds.
$$
\end{lemma}
\vskip 3mm

Please note that the expression
$\int_{s_{0}-\gamma}^{s_{1}+\gamma} s 
\operatorname{vol}_{n}(M_{n}^\circ(s) )\ ds$ is negative.
This lemma means that the volume of $M_{n}^{\circ}$ is concentrated between 
the hyperplanes through $s_{0}e_{n+1}$ and $s_{1}e_{n+1}$.
\par
Although the inequality 
$s_{0}<s_{1}$ follows from the above computations,  it is comforting to verify this
directly. The inequality $s_{0}<s_{1}$ is equivalent to
$$
(a+b)\sqrt{e}(2-\sqrt{\pi})
>0
$$
which holds for all positive $a$ and $b$.
\vskip 3mm

\noindent
{\bf Proof.}
(i)
By Lemma \ref{Kegel2}, 
for all $\gamma>0$ there is $n_{0}$ such that for all $n$
with $n\geq n_{0}$ and all $s$ with
$$
-\frac{1}{a}
\leq s\leq\frac{1-\sqrt{e}}{b+a\sqrt{e}}-\gamma
=s_{0}-\gamma
$$
we have
$$
\frac{1}{2}
(1+s \ a)^{n}\operatorname{vol}_{n}(B_2^n)^{2}
\leq
\operatorname{vol}_{n}(M^\circ(s))
\leq(1+s \ a)^{n}\operatorname{vol}_{n}( B_2^n)^{2}
$$
Therefore
\begin{eqnarray*}
&&\int_{s_{0}-2\gamma}
^{s_{0}-\gamma}  
\operatorname{vol}_{n}(M_{n}^\circ(s)) \ ds   \\
&&\geq \frac{1}{2}\int_{s_{0}-2\gamma}
^{s_{0}-\gamma}  (1+s \ a)^{n}\operatorname{vol}_{n}(B_2^n)^{2} \ ds
\\
&&=\frac{1}{2a(n+1)}| B_2^n|^{2}
\left(\left(1+a\left(s_{0}-\gamma\right)\right)^{n+1}-
\left(1+a\left(s_{0}-2\gamma\right)\right)^{n+1}\right)   \\
&&=\frac{1}{2a(n+1)}| B_2^n|^{2}
\left(1+a\left(s_{0}-\gamma\right)\right)^{n+1}
\left(1-
\left(\frac{1+a\left(s_{0}-2\gamma\right)}
{1+a\left(s_{0}-\gamma\right)}\right)^{n+1}\right) \ . 
\end{eqnarray*}
For sufficiently large $n$
$$
\int_{s_{0}-2\gamma}
^{s_{0}-\gamma}  
\operatorname{vol}_{n}(M_{n}^\circ(s)) \ ds
\geq\frac{1}{4a(n+1)}| B_2^n|^{2}
\left(1+a\left(s_{0}-\gamma\right)\right)^{n+1}.
$$
On the other hand, by Lemma \ref{Kegel2} 
\begin{eqnarray*}
&&\int_{-\frac{1}{a}}^{s_{0}-2\gamma}  
\operatorname{vol}_{n}(M_{n}^\circ(s)) \ ds   \\
&&\leq\int_{-\frac{1}{a}}^{s_{0}-2\gamma}  
(1+s \ a)^{n}\operatorname{vol}_{n}(B_2^n)^{2} \ ds    \\
&&=\frac{1}{a(n+1)}\operatorname{vol}_{n}(B_2^n)^{2}
\left(1+a\left(s_{0}-2\gamma\right)\right)^{n+1}  \\
&&\leq\frac{1}{a(n+1)}\operatorname{vol}_{n}(B_2^n)^{2}
\left(1+a\left(s_{0}-\gamma\right)\right)^{n+1}
\left(1-\frac{a\gamma}
{1+a\left(s_{0}-\gamma\right)}\right)^{n+1} . 
\end{eqnarray*}
Therefore, for sufficiently large $n$
\begin{equation}\label{ConcentrationM1}
\int_{-\frac{1}{a}}^{s_{0}-2\gamma}  
\operatorname{vol}_{n}(M_{n}^\circ(s)) \ ds
\leq 4\left(1-\frac{a\gamma}
{1+a\left(s_{0}-\gamma\right)}\right)^{n+1} 
\int_{s_{0}-2\gamma}
^{s_{0}-\gamma}  
\operatorname{vol}_{n}(M_{n}^\circ(s)) \ ds.
\end{equation}
Thus
\begin{equation*}
\int_{-\frac{1}{a}}^{s_{0}-2\gamma}  
 \operatorname{vol}_{n}(M^\circ(s)) \ ds
\leq
\frac{\gamma}{2}
\int_{s_{0}-2\gamma}
^{s_{0}-\gamma}   \operatorname{vol}_{n}(M^\circ(s)) \ ds.
\end{equation*}
Now we consider the interval $[s_{1},\frac{1}{b}]$.
By Lemma \ref{Kegel2} 
\begin{eqnarray*}
&&\int_{s_{1}+\gamma}^{s_{1}+2\gamma}
\operatorname{vol}_{n}(M_{n}^{\circ}(s))ds   \\
&&\geq\int_{s_{1}+\gamma}^{s_{1}+2\gamma}2^{n-1}(1-s \
b)^{n}\operatorname{vol}_{n}(B_{1}^{n})ds  \\ 
&&=2^{n-1}
\frac{\operatorname{vol}_{n}(B_{1}^{n})}{b(n+1)}\left((1-b(s_{1}+\gamma))^{n+1}-(1-b(s_{1}+2\gamma))^{n+1}
\right)
\\
&&=2^{n-1}\operatorname{vol}_{n}(B_{1}^{n})
\frac{(1-b(s_{1}+\gamma))^{n+1}}{b(n+1)}
\left(1-\left(1-\frac{b\gamma}{1-b(s_{1}+\gamma)}\right)^{n+1}
\right)  
\end{eqnarray*}
Therefore, for sufficiently large $n$
$$
\int_{s_{1}+\gamma}^{s_{1}+2\gamma}
\operatorname{vol}_{n}(M_{n}^{\circ}(s))ds
\geq2^{n-2}\operatorname{vol}_{n}(B_{1}^{n})
\frac{(1-b(s_{1}+\gamma))^{n+1}}{b(n+1)}.
$$
On the other hand, by Lemma \ref{Kegel2} 
\begin{eqnarray*}
&&\int_{s_{1}+2\gamma}^{\frac{1}{b}}
\operatorname{vol}_{n}(M_{n}^{\circ}(s))ds  \\
&&\leq  \int_{s_{1}+2\gamma}^{\frac{1}{b}}2^{n}(1-s \
b)^{n}\operatorname{vol}_{n}(B_{1}^{n})ds   \\
&&=2^{n}\frac{1}{b(n+1)}(1-
b(s_{1}+2\gamma))^{n+1}\operatorname{vol}_{n}(B_{1}^{n})  \\
&&=2^{n}\frac{1}{b(n+1)}
(1-b(s_{1}+\gamma))^{n+1}
\left(1-\frac{b\gamma}{1-b(s_{1}+\gamma)}\right)^{n+1}
\operatorname{vol}_{n}(B_{1}^{n}).
\end{eqnarray*}
Thus
\begin{equation}\label{ConcentrationM2}
\int_{s_{1}+2\gamma}^{\frac{1}{b}}
\operatorname{vol}_{n}(M_{n}^{\circ}(s))ds
\leq 4
\left(1-\frac{b\gamma}{1-b(s_{1}+\gamma)}\right)^{n+1}
\int_{s_{1}+\gamma}^{s_{1}+2\gamma}
\operatorname{vol}_{n}(M_{n}^{\circ}(s))ds\ .
\end{equation}
Since $s_{1}<0$, for sufficiently big $n$
\begin{equation}\label{ConcentrationM3}
\int_{s_{1}+2\gamma}^{\frac{1}{b}}
\operatorname{vol}_{n}(M_{n}^{\circ}(s))ds
\leq \frac{\gamma}{2}
\int_{s_{1}+\gamma}^{s_{1}+2\gamma}
\operatorname{vol}_{n}(M_{n}^{\circ}(s))ds.
\end{equation}
It is left to pass to a new $\gamma$.
\vskip 3mm
(ii) By (\ref{ConcentrationM1})
\begin{eqnarray*}
&&\int_{-\frac{1}{a}}^{s_{0}-2\gamma}  |s|
\operatorname{vol}_{n}(M_{n}^\circ(s)) \ ds   \\
&&\leq 4\max\left\{\frac{1}{a},\frac{1}{b}\right\}
\left(1-\frac{a\gamma}
{1+a\left(s_{0}-\gamma\right)}\right)^{n+1} 
\int_{s_{0}-2\gamma}
^{s_{0}-\gamma}  
\operatorname{vol}_{n}(M_{n}^\circ(s)) \ ds\ .
\end{eqnarray*}
It is left to choose $n$ sufficiently big.
The other estimate is done in the same way using 
(\ref{ConcentrationM3}).
$\square$
\vskip 4mm

\noindent
{\bf Proof of Proposition \ref{thm:distance}.}
(i) is proved in the remark after the statement of 
Proposition \ref{thm:distance}.
We show
(ii). By Definition
$$
g(M_{n}^{\circ})(n+1)
=\frac{\int_{-\frac{1}{a}}^{\frac{1}{b}}
s\operatorname{vol}_{n}(M_{n}^{\circ}(s))ds}
{\int_{-\frac{1}{a}}^{\frac{1}{b}}\operatorname{vol}_{n}(M_{n}^{\circ}(s))ds}.
$$
Therefore, by Lemma \ref{ConcentrationM}
\begin{eqnarray*}
&&\left|
g(M_{n}^{\circ})(n+1)
-\frac{\int_{s_{0}-\gamma}^{s_{1}+\gamma}
s\operatorname{vol}_{n}(M_{n}^{\circ}(s))ds}
{\int_{-\frac{1}{a}}^{\frac{1}{b}}\operatorname{vol}_{n}(M_{n}^{\circ}(s))ds}
\right|   \\
&&\leq
\frac{\int_{-\frac{1}{a}}^{s_{0}-\gamma}
|s|\operatorname{vol}_{n}(M_{n}^{\circ}(s))ds}
{\int_{-\frac{1}{a}}^{\frac{1}{b}}\operatorname{vol}_{n}(M_{n}^{\circ}(s))ds}
+\frac{\int_{s_{1}+\gamma}^{\frac{1}{b}}
|s|\operatorname{vol}_{n}(M_{n}^{\circ}(s))ds}
{\int_{-\frac{1}{a}}^{\frac{1}{b}}\operatorname{vol}_{n}(M_{n}^{\circ}(s))ds}
\leq 2\gamma.
\end{eqnarray*}
Thus
$$
g(M_{n}^{\circ})(n+1)-2\gamma
\leq\frac{\int_{s_{0}-\gamma}^{s_{1}+\gamma}
s\operatorname{vol}_{n}(M_{n}^{\circ}(s))ds}
{\int_{-\frac{1}{a}}^{\frac{1}{b}}\operatorname{vol}_{n}(M_{n}^{\circ}(s))ds}
\leq g(M_{n}^{\circ})(n+1)+2\gamma.
$$
Since $s_{1}<0$,
we may assume that $s_{1}+\gamma<0$. Therefore 
$$
\int_{s_{0}-\gamma}^{s_{1}+\gamma}
s\operatorname{vol}_{n}(M_{n}^{\circ}(s))ds<0
$$
and by Lemma \ref{ConcentrationM} (i)
$$
g(M_{n}^{\circ})(n+1)-2\gamma
\leq\frac{\int_{s_{0}-\gamma}^{s_{1}+\gamma}
s\operatorname{vol}_{n}(M_{n}^{\circ}(s))ds}
{\int_{-\frac{1}{a}}^{\frac{1}{b}}\operatorname{vol}_{n}(M_{n}^{\circ}(s))ds}
\leq\frac{\int_{s_{0}-\gamma}^{s_{1}+\gamma}
s\operatorname{vol}_{n}(M_{n}^{\circ}(s))ds}
{(1+\gamma)\int_{s_{0}-\gamma}^{s_{1}+\gamma}
\operatorname{vol}_{n}(M_{n}^{\circ}(s))ds}
\leq\frac{s_{1}+\gamma}{1+\gamma}\ .
$$
On the other hand
$$
s_{0}-\gamma
\leq
\frac{\int_{s_{0}-\gamma}^{s_{1}+\gamma}
s\operatorname{vol}_{n}(M_{n}^{\circ}(s))ds}
{\int_{s_{0}-\gamma}^{s_{1}+\gamma}\operatorname{vol}_{n}(M_{n}^{\circ}(s))ds}
\leq
\frac{\int_{s_{0}-\gamma}^{s_{1}+\gamma}
s\operatorname{vol}_{n}(M_{n}^{\circ}(s))ds}
{\int_{-\frac{1}{a}}^{\frac{1}{b}}\operatorname{vol}_{n}(M_{n}^{\circ}(s))ds}
\leq g(M_{n}^{\circ})(n+1)+2\gamma\ .
$$
Therefore,
$$
s_{0}-3\gamma
\leq g(M_{n}^{\circ})(n+1)
\leq \frac{s_{1}+\gamma}{1+\gamma}+2\gamma.
$$
We apply now these estimates to the convex body
$$
M_n
=\operatorname{co}\left[(K,-1),\left(L,\frac{1}{e-1}\right)\right].
$$
The centroid of $M_{n}$ is $(0,\delta_{n})$ with 
$\lim_{n\to\infty}\delta_{n}=0$. We get
\begin{eqnarray*}
&&M_{n}^{(0,\delta_{n})}
=\\
&&\left\{(z,s)\left|\forall x\in K:\langle z,x\rangle-s(1+\delta_{n})\leq1
\hskip 1mm\mbox{and}\hskip 1mm
\forall y\in L:\langle z,y\rangle+s\left(\frac{1}{e-1}-\delta_{n}\right)\leq1)
\right.\right\}.
\end{eqnarray*}
It is left to apply the above estimates to $M_{n}^{\circ}$ with
$a=1+\delta_{n}$ and $b=\frac{1}{e-1}-\delta_{n}$.
$\square$

\vskip 2mm 

Mathieu Meyer \\
  {\small        Universit\'{e} de Paris Est - Marne-la-Vall\'{e}e}\\
    {\small      Equipe d'Analyse et de Mathematiques Appliqu\'{e}es}\\
   {\small       CitŽ Descartes - 5, bd Descartes }\\
    {\small      Champs-sur-Marne
          77454 Marne-la-Vall\'{e}e,  France} \\
          {\small \tt mathieu.meyer@univ-mlv.fr} \\

   Carsten Sch\"utt \\
     {\small       Christian Albrechts Universit\"at }\\
        {\small    Mathematisches Seminar }\\
        {\small    24098 Kiel, Germany} \\
         {\small \tt schuett@math.uni-kiel.de}   \\

     \and Elisabeth Werner\\
{\small Department of Mathematics \ \ \ \ \ \ \ \ \ \ \ \ \ \ \ \ \ \ \ Universit\'{e} de Lille 1}\\
{\small Case Western Reserve University \ \ \ \ \ \ \ \ \ \ \ \ \ UFR de Math\'{e}matique }\\
{\small Cleveland, Ohio 44106, U. S. A. \ \ \ \ \ \ \ \ \ \ \ \ \ \ \ 59655 Villeneuve d'Ascq, France}\\
{\small \tt elisabeth.werner@case.edu}\\ \\

\end{document}